\def\Q{{\mathbb Q}}
\def\Z{{\mathbb Z}}

\def\C{{\mathbb C}}

\def\Gal{\mathrm{Gal}}

\def\Aut{\mathrm{Aut}}
\def\Hom{\mathrm{Hom}}
\def\Sin{\mathrm{Sin}}
\def\Norm{\mathrm{Norm}}
\def\Cent{\mathrm{Cent}}
\def\Out{\mathrm{Out}}
\def\Inn{\mathrm{Inn}}  
\def\Diff{\mathrm{Diff}}
\def\Spec{\mathrm{Spec}}
\def\pro{\mathrm{pro}}
\def\I{{I}}

\def\s{{\mathcal S}}
\def\t{{\mathcal T}}

\def\sk{\mathrm{sk}}
\def\cosk{\mathrm{cosk}}

\def\dim{\mathrm{dim}}
\def\O{{\mathcal O}}

\def\hf{{\mathfrak S}}
\def\vf{{\mathfrak h}}

\def\T{{\mathcal T}}
\def\E{{\mathfrak E}}
\def\m{{\mathcal M}}
\def\h{{\mathcal H}}
\def\O{{\mathcal O}}
\def\I{{\mathcal I}}
\def\r{{\mathcal R}}
\def\e{{\}_{et}}}

\newcommand{\ms}{\medskip}

\documentstyle[amssymb,amsfonts,amscd,11pt]{amsart}
\input xypic

\newtheorem{theorem}{Theorem}[section]
\newtheorem*{itheo}{Theorem}

\newtheorem{proposition}[theorem]{Proposition}
\newtheorem{corollary}[theorem]{Corollary}

\theoremstyle{definition}
\newtheorem{definition}[theorem]{Definition}

\newtheorem{example}[theorem]{Example}

\theoremstyle{remark}

\numberwithin{equation}{section}

\begin{document}

\title[Etale Homotopy of Moduli of Curves with Symmetries]
{Etale Homotopy Types of Moduli Stacks\\
of Algebraic Curves with Symmetries\\}

\author[Paola Frediani]
{Paola Frediani}
\address{Dipartimento di Matematica\\
Universit\`a degli Studi di Pavia\\
Via Ferrata 1, I-27100 Pavia, Italy}
\email{frediani@@dimat.unipv.it}

\author[Frank Neumann]
{Frank Neumann}
\address{Department of Mathematics and Computer Science\\
University of Leicester\\
University Road, Leicester LE1 7RH, England, UK}
\email{fn8@@mcs.le.ac.uk}




\begin{abstract}
{Using the machinery of etale homotopy theory \`{a} la Artin-Mazur we determine the etale homotopy types of 
moduli stacks over $\bar{\Q}$ parametrizing families of algebraic curves of genus $g \geq 2$ endowed with an action 
of a finite group $G$ of automorphisms, which comes with a fixed embedding in the mapping class group $\Gamma_g$, 
such that in the associated complex analytic situation the action of $G$ is precisely the differentiable action
induced by this specified embedding of $G$ in $\Gamma_g$.}
\end{abstract}

\maketitle

\hfill {\it Dedicated to Professor Hyman Bass}\\
\hbox to .99\textwidth{\hfill {\it on the occasion of his 70th birthday.}}

\section*{Introduction}

\ms\noindent Let $\m_g$ be the moduli stack of families of proper smooth algebraic curves of genus $g$ with
$g\geq 2$. In a fundamental paper Oda \cite{O} determined the etale homotopy type of $\m_g \otimes \bar{\Q}$,
the moduli stack representing the restriction of the moduli stack to the subcategory of schemes over
$\bar{\Q}$. The etale homotopy type a la Artin-Mazur of the moduli stack
$\m_g \otimes \bar{\Q}$ is given as the profinite Artin-Mazur completion of the Eilenberg-MacLane
space $K(\Gamma_g, 1)$, where $\Gamma_g$ is the Teichm\"uller modular or mapping class group of compact Riemann
surfaces of genus $g$. In fact, Oda proved this result for all moduli stacks $\m_{g,n}$ of families of proper
smooth algebraic curves of genus $g$ with $n$ distinct ordered points and $2g+ n > 2$. In this more general
context, the group $\Gamma_g$ is replaced with the Teichm\"uller modular group $\Gamma_{g, n}$ of compact
Riemann surfaces of genus $g$ with $n$ punctures. These results of Oda on the etale fundamental group
of the moduli stack of algebraic curves are of great importance for the geometry and arithmetic of moduli spaces of 
algebraic curves and for the study of geometric Galois actions as they allow to relate the absolute
Galois group $\Gal(\bar{\Q}/\Q)$ to the profinite completion $\Gamma_{g,n}^{\wedge}$ of the mapping class groups 
$\Gamma_{g,n}$ which is important for realizing Grothendieck's ideas on a ``Lego-Teichm\"uller game'' \cite{G1}. 
We refer to the excellent articles in the book edited by Schneps and Lochak \cite{SL} and the article by 
Matsumoto \cite{Ma} for good overviews on this fascinating cycle of ideas.  

\ms\noindent In this article we embark in a different direction. We study the etale homotopy types of certain moduli stacks $\m_{g, G}$ 
parametrizing families of proper smooth algebraic curves of genus $g\geq 2$ with certain prescribed fixed finite subgroups of automorphisms.
It turns out that a similar theorem like that of Oda for $\m_g$ holds with some modifications also in this context. We outline 
very briefly the argument here. We will study the following general situation: Let $G$ be a finite group with a fixed embedding 
in the mapping class group $\Gamma_g$ with $g\geq 2$. Let now $\m_{g, G}$ be the moduli stack of those families of proper smooth
algebraic curves of genus $g\geq 2$ which have $G$ as a subgroup of its automorphism groups and which are ``keeping track of the
differentiable action of $G$ as a subgroup of $\Gamma_g$ in the associated complex analytic situation''.
More precisely, the moduli stack $\m_{g, G}$ of algebraic curves with symmetries is defined as the category fibred in groupoids over 
the category of schemes over $\Q$ defined by its category of sections as follows: For a scheme $S$ over $\Q$, the objects of $\m_{g, G}(S)$ 
are families of algebraic curves over $S$ with an action of the finite group $G$ such that the differentiable action of $G$ on the associated 
complex analytic family of Riemann surfaces of genus $g$ is precisely the one induced by the 
specified embedding of $G$ in $\Gamma_g$. The morphisms of $\m_{g, G}(S)$ are those morphisms of families of algebraic curves over $S$ which are 
compatible with these data. It turns out that $\m_{g, G}$ is in fact, like $\m_g$, an algebraic stack in the sense of Deligne and Mumford \cite{DM}.

\ms\noindent In order to determine the etale homotopy type $\{\m_{g, G}\otimes \bar{\Q}\e$ of the stack 
$\m_{g, G}\otimes \bar{\Q}$ we compare the algebraic with the associated complex analytic situation over
the complex numbers. There we determine using Teichm\"uller theory the analytic or classical topological homotopy type of the
complex analytification $\m_{g, G}^{an}$ of the stack $\m_{g, G}$ parametrizing families of compact Riemann surfaces
of genus $g\geq 2$ with differentiable action of $G$ as given by the fixed embedding of $G$ in $\Gamma_g$. 
Moduli spaces of Riemann surfaces with symmetries and their related Teichm\"uller theory were studied systematically from a 
complex geometry point of view by Gonzal\'ez-D\'{\i}ez and Harvey \cite{GDH}. They studied irreducible subvarieties of the moduli space of Riemann
surfaces of genus $g\geq 2$ and their normalizations characterized by specifying a finite subgroup $G$ of the mapping class group whose action
on the surfaces is fixed geometrically. From their results it follows that in our notation, the complex analytification
$\m_{g, G}^{an}$ of the stack $\m_{g, G}$ is basically the orbifold given as the quotient $T_g^G/\Gamma_{g,G}$, where 
$T_g^G$ is the fixed point locus of the classical Teichm\"uller space $T_g$, representing the moduli functor for 
families of marked compact Riemann surfaces of genus $g$ and where the modular group $\Gamma_{g,G}$ is the normalizer
of the specified finite group $G$ in the mapping class group $\Gamma_g$. Especially, it follows that the stack
$\m_{g, G}^{an}$ is actually a normal complex analytic space, the normalization of a substack ${\mathcal M'}_{g, G}^{an}$ of 
the analytification $\m_g^{an}$ of the moduli stack $\m_g$.

\ms\noindent The solution of the Nielsen realization problem by Kerkhoff \cite{K} and others (cf. also the articles by Tromba \cite{Tr} and Catanese\cite{Ca1}) 
implies that the Teichm\"uller space $T_g^G$ is in fact a contractible topological space, which allows to determine
the classical homotopy type of the stack $\m_{g, G}^{an}$ and using a general comparison theorem comparing the etale
and classical homotopy types over the complex numbers essentially due to Artin-Mazur \cite{AM} (cf. also Cox \cite{C} and Friedlander \cite{F}) 
we finally derive our main theorem.

\ms
\begin{itheo}
There is a weak homotopy equivalence of pro-simplicial sets
$$\{{\mathcal M}_{g, G}\otimes\bar{\Q}\}^{\wedge}_{et}\simeq
 K(\Gamma_{g,G}, 1)^{\wedge}.$$
where $\Gamma_{g, G}$ is the normalizer of the group $G$ in the mapping class group $\Gamma_g$ and $^{\wedge}$ denotes
Artin-Mazur profinite completion.
\end{itheo}

\ms\noindent An interesting special case is given by the moduli stack $\h_g$ of families of hyperelliptic curves, i.e.,
families of algebraic curves with a hyperelliptic involution, which corresponds to the moduli stack $\m_{g, G}$ where
$G$ is simply the group $\Z/2$. This particular case was already mentioned by Oda \cite{O}.

\ms\noindent The paper is organized as follows: In the first chapter we define homotopy types of Deligne-Mumford
stacks using the machinery of Artin-Mazur \cite{AM}. We also compare the etale and classical homotopy types in
the algebraic and complex analytic context. In the second chapter we collect the necessary facts from Teichm\"uller theory for
families of Riemann surfaces with symmetries, introduce the analytic stacks $\m_{g, G}^{an}$ and determine their
classical homotopy types. Finally, in the third chapter we introduce the algebraic stacks $\m_{g, G}$ of families
of algebraic curves with symmetries and determine their etale homotopy types using the results from the second and first chapter.

\ms\noindent
{\it Acknowledgements.} The first author likes to thank all the people 
at the  Centre de Recerca Matem\`atica, Bellaterra, Spain and the members of the Barcelona Algebraic Topology Group for much 
support and a wonderful stay, while holding a post-doc research grant as part of the EU research 
training network Modern Homotopy Theory.\\
\noindent Both authors like to thank Fabrizio Catanese for his interest and support and the people of the former Algebraic Geometry
Seminar in G\"ottingen, where discussion on this work was initiated.\\
\noindent We also like to thank the referee for some useful comments and suggestions.

\section{Homotopy Types of Deligne-Mumford stacks}

\subsection{Homotopy types for locally connected topoi.}

We briefly recall the construction of homotopy types for locally connected topoi following
Artin-Mazur \cite{AM} and Cox \cite{C}. We use the language of topoi here following Moerdijk \cite{Mo} 
and Zoonekynd \cite{Z}. 

\ms\noindent
Let $\E$ be a topos and let $\emptyset$ be its initial object. An object $X$ of $\E$ is {\it connected},
if whenever $X=X_1\coprod X_2$, either $X_1$ or $X_2$ is the initial object $\emptyset$. 
The topos $\E$ is {\it locally connected} if every object of $\E$ is a coproduct of connected objects.

\ms\noindent
Let $\Delta$ be the {\it category of simplices} i.e. the category whose objects are sets 
$[n]=\{0, 1, 2, \ldots n\}$ and whose morphisms are non-decreasing maps. Let also $\Delta_n$ be the
full subcategory of $\Delta$ of all sets $[k]$ with $k\leq n$.

\ms\noindent
Let $\E$ be a locally connected topos. Define $\pi$ to be the {\it connected component functor}
$$\pi: \E \rightarrow (Sets)$$
associating to any object $X$ of $\E$ its set $\pi(X)$ of connected components.

\ms\noindent
A {\it simplicial object} in $\E$ is a functor $X_{\bullet}: \Delta^{op} \rightarrow \E$. The category of
all simplicial objects in the topos $\E$ will be denoted by $\Delta^{op} \E$.

\ms\noindent
The {\it restriction} or {\it $n$-truncation functor} 
$$(-)^{(n)}: \Delta^{op} \E \rightarrow \Delta_n^{op} \E$$
has left and right adjoint functors, the $n${\it -th skeleton} and $n${\it -th coskeleton functor}
$$\sk_n, \; \cosk_n: \Delta_n^{op} \E \rightarrow \Delta^{op} \E.$$

\ms\noindent
A {\it hypercovering} of the topos $\E$ is a simplicial object $U_{\bullet}$ such that the
morphisms
$$U_0 \rightarrow *$$
$$U_{n+1}\rightarrow \cosk_n(U_{\bullet})_{n+1}$$
are epic, where $*$ is the final object of $\E$.

\ms\noindent
If $S$ is a set and $X$ an object of the topos $\E$ define $T\otimes X:= \coprod_{s\in S}X$.
If $S_{\bullet}$ is a simplicial set and $X_{\bullet}$ a simplicial object of $\E$ define 
$$S_{\bullet}\otimes X_{\bullet}: \Delta^{op} \rightarrow \E$$
to be the simplicial object given by $(S_{\bullet}\otimes X_{\bullet})_n =S_n \otimes X_n.$

\ms\noindent
Let $\Delta[m] =\Hom_{\Delta}(-, [m])$ be the standard $m$-simplicial set i.e., the functor $\Delta[m]: \Delta^{op}
\rightarrow (Sets)$ represented by the set $[m]$.

\ms\noindent
Two morphisms $f,g: X_{\bullet}\rightarrow Y_{\bullet}$ are {\it strictly homotopic} if there is a 
{\it strict homotopy} $H: X_{\bullet} \otimes \Delta[1]\rightarrow Y_{\bullet}$ such that the following diagram
is commutative
$$\diagram
X_{\bullet}=X_{\bullet}\otimes\Delta[0] \dto_{id\otimes d^0} \drrto^{f} &   &  \\
X_{\bullet}\otimes\Delta[1] \rrto^{H}                                        &   & Y_{\bullet} \\
X_{\bullet}=X_{\bullet}\otimes\Delta[0] \uto^{id\otimes d^1} \urrto^{g} &   &   \enddiagram$$

\ms\noindent
Two morphisms $f,g: X_{\bullet}\rightarrow Y_{\bullet}$ are {\it homotopic}, if they can be related by a chain
of strict homotopies. {\it Homotopy} is the equivalence relation generated by strict homotopy.

\ms\noindent
Let $\h\r (\E)$ be the {\it homotopy category of hypercoverings of} $\E$, i.e., the category of
hypercoverings of $\E$ and their morphisms up to homotopy.
It turns out that the opposite category $\h\r (\E)^{op}$ is actually a filtering category. The proof
of \cite{AM}, Corollary 8.13 applies verbatim.

\ms\noindent
Let also $\h=\h(\Delta^{op}(Sets))$ be the {\it homotopy category of simplicial sets}.
This category is actually equivalent to the homotopy category of CW-complexes (cf. \cite{GJ}, Chap. I or \cite{BK}, VIII.3). 

\ms\noindent
Further let $\pro-\h$ be the {\it category of pro-objects} in the category $\h$, i.e., the category of (contravariant) functors
$X: \I \rightarrow \h$ from some filtering index category $\I$ to $\h$. We will write normally $X=\{X_i\}_{i\in \I}$ to
indicate that we think of pro-objects $X$ as inverse systems of objects of $\h$ (cf. \cite{AM}, Appendix).

\ms\noindent
Let $\E$ be a locally connected topos. The connected component functor $\pi: \E \rightarrow (Sets)$
induces a functor
$$\Delta^{op}\pi: \Delta^{op}\E \rightarrow \Delta^{op}(Sets).$$
\noindent Passing to homotopy categeories and restricting to hypercoverings of $\E$ gives a functor
$$\pi: \h\r(\E) \rightarrow \pro-\h.$$

\ms\noindent Now we can define the homotopy type of a locally connected topos.

\begin{definition}
Let $\E$ be a locally connected topos. The {\it homotopy type} of $\E$ is given as the following pro-object
in the homotopy category of simplicial sets:
$$\{\E\}= \{\pi(U_{\bullet})\}_{U_{\bullet}\in \h\r(\E)}.$$
\end{definition}

\ms\noindent
This construction is actually functorial with respect to morphisms of topoi, the associated functor $\{-\}$ is also called
the {\it Verdier functor}.

\ms\noindent If $\E$ is a locally connected topos and $x$ a point of the topos $\E$, i.e., a morphism of topoi 
$x: \mathfrak{S}\mathfrak{e}\mathfrak{t} \rightarrow \E$,
one can also define {\it homotopy groups} $\pi_n(\E, x)$ for $n\geq 1$ following Artin-Mazur \cite{AM} 
(cf. also \cite{Mo} and \cite{Z}). In general these homotopy groups $\pi_n(\E, x)$ turn out to be 
pro-groups.

\subsection{Etale homotopy types of Deligne-Mumford stacks.}

\noindent We will now define the etale homotopy type for the topos of sheaves on the small etale site 
of a Deligne-Mumford stack. 
For the definitions and properties of stacks in general we refer to the book of Laumon and Moret-Bailly \cite{LMB}, to the 
appendix of the article of Vistoli \cite{Vi} and to the original article of Deligne-Mumford \cite{DM}. 
For a nice overview we also recommend the article of Fantechi \cite{fa} and on algebraic stacks in particular 
the expository article of G\'omez \cite{Go}.

\ms\noindent 
Let us first recall that an {\it algebraic Deligne-Mumford stack} $\s$ is basically a contravariant $2$-functor
from the big etale site of the category of schemes to the $2$-category of groupoids being a sheaf with effective descent
data such that the diagonal $\Delta: \s\rightarrow \s\times \s$ is
quasi-compact, separated and representable by schemes and such that there
exists a surjective etale $1$-morphism $x: X\rightarrow \s$ where $X$ is a scheme (called {\it atlas}).
We refer to \cite{Go}, Chap. 2 for the exact definitions using the language of $2$-functors.
Equivalently, following the original definition of Deligne-Mumford \cite{DM}, \S4 an algebraic Deligne-Mumford stack is
a category fibred in groupoids over the big etale site of the category of schemes such that isomorphisms are a sheaf,
all descent data are effective, the diagonal is quasi-compact, separated and representable by schemes and there
exists a surjective etale morphism from a scheme.
From a slightly different viewpoint a Deligne-Mumford stack can also be viewed equivalently as an etale groupoid
in the category of schemes, i.e., a groupoid 
$R \begin{array}{c}\rightarrow\\ \rightarrow \end{array} U$ 
in the category $(Sch)$ of schemes such that the two structure morphisms $R \begin{array}{c}\rightarrow\\ \rightarrow \end{array} U$
are etale morphisms (cf. \cite{Vi}, Appendix).

\ms\noindent
Let $\s$ be an algebraic Deligne-Mumford stack. Then we can consider
its small etale site $\s _{et}$. The objects are etale $1$-morphisms $x: X\rightarrow \s$,
where $X$ is a scheme, morphisms are morphisms over $\s$, i.e., commutative
diagrams of $1$-morphisms of the form

$$\diagram
   X \rrto \drto&               & Y\dlto \\
                 & \s &                 \enddiagram $$

\noindent and the coverings are the etale coverings of the schemes. 

\ms\noindent
An important property of algebraic Deligne-Mumford stacks is the existence of
fibre products $X\times_{\s} Y$ in $\s_{et}$ \cite{DM}.

\ms\noindent
We let $\E_{et}=\hf\vf(\s_{et})$ be the associated etale topos, i.e., the category of 
sheaves on the small etale site of the algebraic stack $\s$. It turns out that this topos $\E_{et}$ is
actually a locally connected topos (cf. \cite{Z}, 3.1.)

\ms\noindent 
Now we define the etale homotopy type of an algebraic Deligne-Mumford stack as the Artin-Mazur homotopy
type of its etale topos via the procedure of the previous section 1.1.

\begin{definition}
Let $\s$ be an algebraic Deligne-Mumford stack. The {\it etale homotopy type} $\{\s\}_{et}$ of $\s$ is the
pro-object in the homotopy category of simplicial sets given as
$$\{\s\}_{et}=\{\E_{et}\}.$$
\end{definition}

\ms\noindent The Artin-Mazur homotopy groups $\pi_n(\E_{et}, x)$ will be also denoted as $\pi_n^{et}(\s, x)$ and called
the {\it etale homotopy groups} of the stack $\s$. For $n=1$ this gives the etale fundamental group $\pi_1^{et}(\s, x)$
of $\s$, which is discussed in detail also in Zoonekynd $\cite{Z}$. 

\ms\noindent One should not confuse this Artin-Mazur etale homotopy type with the etale topological type of Friedlander \cite{F}, 4 
or with the rigid etale homotopy type of Hoobler and Rector \cite{HR} where rigid hypercoverings are used instead in the definition.
Friedlander uses in \cite{F} the subscript 'et' for the etale topological type and the subscript 'ht' for the Artin-Mazur etale 
homotopy type. We use 'et' here for the Artin-Mazur etale type in order to be coherent with the notation used in Oda \cite{O}.

\ms\noindent
The following generalization in the context of algebraic stacks of the homotopy descent
theorem for simplicial schemes of Cox \cite{C}, Theorem IV.2 (cf. also \cite{O}, Theorem 3), is straightforward and allows to determine the 
etale homotopy type of a Deligne-Mumford stack directly from that of a hypercovering $X_{\bullet}$, i.e., from a simplicial scheme over 
$\s$ which can often be constructed directly.

\begin{theorem}
If $\s$ is an algebraic Deligne-Mumford stack with etale site $\s_{et}$ and
$X_{\bullet}$ a simplicial scheme which is a hypercovering of $\E_{et}$,
then there is a weak homotopy equivalence of pro-simplicial sets
$$\{\s\e \simeq \{X_{\bullet}\e.$$
where $\{X_{\bullet}\e$ is the etale homotopy type of the simplicial scheme $X_{\bullet}$.
\end{theorem}

\ms\noindent 
The etale homotopy type $\{X_{\bullet}\e$ of the simplicial scheme $X_{\bullet}$ over $\s$ is defined here as in Cox \cite{C},
Chap. III, where its relation with the etale topological type of Friedlander \cite{F}, 4 is also explained. It turns out that if
the simplicial scheme $X_{\bullet}$ which is a hypercovering of $\E_{et}$ has the property that each $X_n$ is quasi-projective
then the etale Artin-Mazur homotopy type and the etale topological type are weakly equivalent in $\pro-\h$. 

\ms\noindent
The existence of such a hypercovering follows basically from the existence of a surjective, etale morphism 
$x: X\rightarrow \s$ where $X$ is a scheme, which is intrinsic in the definition of a Deligne-Mumford stack $\s$. 
The hypercovering can then be constructed via itered fiber products along this etale covering morphism.

\subsection{Analytic Homotopy Types of Deligne-Mumford stacks.}

\ms\noindent
We will consider now Deligne-Mumford stacks in the context of complex analytic spaces. An excellent overview
on analytic stacks in general and analytification as well as GAGA type theorems can be found in Toen \cite{T}, 
Chap. 5. and we refer to that article for exact definitions and properties of analytic stacks and analytification.

\ms\noindent In analogy, an {\it analytic Deligne-Mumford stack} is basically a contravariant $2$-functor from the big site of local isomorphisms 
of the category of complex analytic spaces to the $2$-category of groupoids being a sheaf with effective descent data such that 
the diagonal $\Delta: \t\rightarrow \t\times \t$ is
finite and representable by complex analytic spaces and such that there
exists a surjective etale $1$-morphism $x: X\rightarrow \t$ where $X$ is here
a complex analytic space. Again, we could think of an analytic Deligne-Mumford stack as an etale
groupoid in the category $(AnSp)$ of complex analytic spaces or we could use also here the language of
categories fibred in groupoids.

\ms\noindent
For such an analytic Deligne-Mumford stack $\t$ we can define its site of local isomorphisms
$\t_{cl}$ given by local isomorphisms $x: X\rightarrow \t$, where $X$ is
a complex analytic space and morphisms are morphisms of analytic spaces over
$\t$ and the coverings are given by families of local isomorphisms (cf. Mumford \cite{M},
and Cox \cite{C}, Chap. IV, \S3). 

\ms\noindent
We let $\E_{cl}=\hf\vf(\t_{cl})$ be the associated topos of local isomorphisms, i.e., the category of 
sheaves on the small site of local isomorphisms of the analytic stack $\t$. It is again locally connected
and we can define the analytic homotopy type of $\t$.

\begin{definition}
Let $\t$ be an analytic Deligne-Mumford stack. The {\it analytic} or {\it classical homotopy type} $\{\t\}_{cl}$ of $\t$ is the
pro-object in the homotopy category of simplicial sets given as
$$\{\t\}_{cl}=\{\E_{cl}\}.$$
\end{definition}

\ms\noindent
For an analytic Deligne-Mumford stack $\t$ we have actually a very explicit description
of the classical homotopy type $\{\t\}_{cl}$ using hypercoverings and the underlying topological spaces. 
Namely, if $\t$ is an analytic Deligne-Mumford stack and $X_{\bullet}$ is a simplicial analytic
space, which is also a hypercovering of the topos $\E_{cl}$, then similarly as in the case af algebraic Deligne-Mumford
stacks there is a weak homotopy equivalence 
$$\{\t\}_{cl} \simeq \{X_{\bullet}\}_{cl}$$
\noindent where the classical homotopy type $\{X_{\bullet}\}_{cl}$ of the small site of local isomorphisms of the hypercovering $X_{\bullet}$
is given as  $\{X_{\bullet}\}_{cl} = \Delta \Sin (X_{\bullet})$. Here $\Sin (X_{\bullet})$ is the bisimplicial set 
given in bidegree $s,t$ by $\Sin_t (X_s)$ and $\Delta$ is the diagonal functor (cf. \cite{C}, Chap. IV,  \S3 and \cite{F}, chap. 8).

\ms\noindent
Using the canonical homotopy equivalence
$$\Delta \Sin (X_{\bullet})\stackrel{\simeq}\longrightarrow \Sin (|X_{\bullet}|)$$
where $|X_{\bullet}|$ denotes here the geometric realization of the simplicial space $X_{\bullet}$, we have that
$$\{\t\}_{cl} \simeq \Sin (|X_{\bullet}|).$$
\noindent Therefore we could have defined the classical homotopy type of an analytic
Deligne-Mumford stack directly using a complex analytic hypercovering. In \cite{O} it
is shown that this is independent of the actual choice of the hypercovering $X_{\bullet}$. 

\ms\noindent
Let now $\bar{\Q}\hookrightarrow \C$ be an embedding of the
algebraic closure of the rationals into the complex numbers, then
for any algebraic Deligne-Mumford stack $\s$ over $\bar{\Q}$
let $\s^{an}$ be the associated analytic Deligne-Mumford stack (cf. \cite{T}, Chapitre 5).
Similarly for any scheme $X$ over $\bar{\Q}$, let $X^{an}$ denote the complex
analytic space associated with the $\C$-valued points $X(\C)$ of $X$.

\ms\noindent
As the etale and analytic homotopy types are determined by a hypercovering of their respective topoi, 
it is important to have a good comparison theorem. We have the following general comparison theorem 
between etale and classical homotopy types of Deligne-Mumford stacks. We refer to Artin-Mazur \cite{AM}, \S3
and Friedlander \cite{F}, Chap. 6 for the Artin-Mazur profinite completion functor and its properties.

\begin{theorem}
Let $\s$ be an algebraic Deligne-Mumford stack over $\bar{\Q}$ and
$X_{\bullet}$ a simplicial scheme which is of finite type over $\bar{\Q}$ and a hypercovering
of the topos $\E_{et}=\hf\vf(\s_{et})$. If $X^{an}_{\bullet}$ is the associated simplicial
complex analytic space of $X_{\bullet}$, then there is weak homotopy equivalence of
pro-simplicial sets
$$\{\s\}^{\wedge}_{et}\simeq \Sin (|X^{an}_{\bullet}|)^{\wedge}.$$
\noindent where $^{\wedge}$ denotes the Artin-Mazur profinite completion functor on the homotopy category
of simplicial sets.
\end{theorem}

\ms\noindent {\bf Proof.} This follows directly from the above considerations and the general comparison theorem for 
homotopy types of simplicial schemes (cf. Cox \cite{C}, Theorem IV.8 or Friedlander \cite{F}, Theorem 8.4), because 
the etale homotopy type of a Deligne-Mumford stack is completely determined by a hypercovering, i.e., by a simplicial scheme. \qed

\section{Moduli stacks of families of Riemann surfaces with symmetries and their homotopy types}

\subsection{Recollections on classical Teichm\"uller theory}

\ms
We will first recall some notions and fundamental facts from the Teichm\"uller theory of families of Riemann surfaces 
(cf. also \cite{O}).

\ms\noindent
Let $C_g$ be a compact connected Riemann surface of genus $g$ and $\Pi_g$ be the abstract group with generators
$$\alpha_1, ..., \alpha_g, \beta_1,...,\beta_g$$
and relation $[\alpha_1, \beta_1]...[\alpha_g,\beta_g]=1$. The abstract group $\Pi_g$ is isomorphic to the fundamental group of the
Riemann surface $C_g$.
Recall that for the group $\Out(\Pi_g)$ of outer automorphisms of $\Pi_g$ we have $\Out(\Pi_g) = \pi_0(\Diff(\Sigma_g))$, where $\Sigma_g$ 
is a differentiable model for the Riemann surface of genus $g$. 

\ms\noindent
The group $\pi_0(\Diff(\Sigma_g))$ acts properly discontinuously on the Teichm\"uller space ${T}_g = M(\Sigma_g)/\Diff_0(\Sigma_g)$, where 
$M(\Sigma_g)$ is the moduli space of complex analytic structures on $\Sigma_g$ and $\Diff_0(\Sigma_g)$ the topological group of all diffeomorphisms 
which are isotopic to the identity.
The action of $\pi_0(\Diff(\Sigma_g))$ on $M(\Sigma_g)$ is defined in the following way. Given an element $f$ of $\Diff(\Sigma_g)$ the induced
element on $M(\Sigma_g)$ is given by $f^*: M(\Sigma_g) \rightarrow  M(\Sigma_g)$, where for a complex analytic structure $X \in M(\Sigma_g)$, 
$f^*(X)$ is that complex analytic structure of $\Sigma_g$ for which the map $f: (\Sigma_g, f^*(X)) \rightarrow (\Sigma_g, X)$ is either holomorphic or 
antiholomorphic, depending on whether $f$ is orientation preserving or not.

\ms\noindent
We denote by $\Diff^+(\Sigma_g)$ the topological group of all orientation preserving diffeomorphisms of $\Sigma_g$ and 
by $\Gamma_g = \pi_0(\Diff^+(\Sigma_g))$ the mapping class group of genus $g$.

\ms\noindent
Let now $C_g$ be a compact connected Riemann surface of genus $g \geq 2$ and $c_0$ a point in $C_g$. Algebraically $H^2(\Pi_g, \Z) \cong \Z$ and 
we choose a generator $e$ for $H^2(\Pi_g, \Z)$. 

\ms\noindent
We consider now any isomorphism $\phi: \pi_1(C_g,c_0) \stackrel{\cong} \rightarrow \Pi_g$ such that the induced isomorphism of cohomology groups 
$$H^2(C_g, \Z) = H^2(\pi_1(C_g, c_0), \Z) \cong H^2(\Pi_g, \Z)$$
maps the orientation class of $C_g$ to the given generator $e$ of $H^2(\Pi_g, \Z)$.
We say that two such isomorphisms $\phi_1$ and $\phi_2$ are equivalent if and only if there exists an inner automorphism $\theta$ of $\Pi_g$ 
such that $\phi_2 = \theta \circ \phi_1$. An equivalence class $[\phi]$ of such isomorphisms is called a {\it marking} of the Riemann surface 
$C_g$ and the pair $(C_g, [\phi])$ a {\it marked} Riemann surface of genus $g$.

\ms\noindent
In order to study the moduli problem for families of marked Riemann surfaces, we recall the following basic definitions. In what follows
we will always assume that $g\geq 2$.

\ms
\begin{definition}
Let $S$ be a complex analytic space. A {\it family of complex analytic curves of genus $g$ over $S$} is a proper, smooth surjective holomorphic 
morphism $p: C\rightarrow S$ between complex analytic spaces such that for every point $s$ in $S$ the fibre $C_s=p^{-1}(s)$ 
is a compact connected Riemann surface of genus $g$.
\end{definition}

\ms\noindent
The family version of a marked Riemann surfaces is given by introducing locally constant families of markings.

\ms
\begin{definition}
Let $p:C\rightarrow S$ be a family of complex analytic curves of genus $g$ over $S$. A {\it locally constant family of markings}
of $p$ is an isomorphism up to conjugation of local systems of groups
$$\phi: \Pi_1(C/S)\cong \Pi_g\times S$$
where $\Pi_1(C/S)$ is the local system of the fundamental groups $\pi_1(C_s, c_0)$ of the fibers $C_s$ and $\Pi_g\times S$ is the constant
local system of groups over $S$ with fibers $\Pi_g$ such that $\phi$ induces an isomorphism $R^2p_*\Z\cong H^2(\Pi_g, \Z)\times S$ compatible
with the orientations at each fiber $C_s$. A {\it family of marked complex analytic curves of genus $g$} is a pair $(p:C\rightarrow S, \phi)$, where $p:C\rightarrow S$ is a
family of complex analytic curves over a complex analytic space $S$ and $\phi$ a locally constant family of markings of $p$.
\end{definition}

\ms
\begin{definition}
Let $(p:C\rightarrow S, \phi)$ and $(p':C'\rightarrow S, \phi')$ be two families of marked complex analytic curves
of genus $g$ over $S$. An {\it isomorphism} between $(p:C\rightarrow S, \phi)$ and $(p':C'\rightarrow S, \phi')$ is an $S$-isomorphism $\alpha: C\rightarrow C'$ which is compatible with the families $\phi$ and $\phi'$ of markings.
\end{definition}

\ms\noindent
We can now describe the moduli problem for families of marked Riemann surfaces. 

\ms
\begin{definition}
Let $\t_g: (AnSp)\rightarrow (Sets)$ be the contravariant functor from the category of complex analytic spaces to
the category of sets, which associates to every complex analytic space
$S$ the set of isomorphism classes of families $(p:C\rightarrow S, \phi)$ of marked complex analytic curves of genus $g$ over $S$
and to every holomorphic morphism $f: S'\rightarrow S$ the map between isomorphism classes induced by the base change with $f$.
\end{definition} 

\ms\noindent
From classical Teichm\"uller theory we know that this functor is representable, i.e., the moduli problem has a fine solution (cf. \cite{G}).

\ms
\begin{theorem}
The moduli problem for families of marked Riemann surfaces has a fine solution, i.e., the moduli functor $\t_g$ is representable by a 
complex analytic space $T_g$, called the {\it Teichm\"uller space}.
\end{theorem} 

\subsection{Teichm\"uller theory for families of marked Riemann surfaces with symmetries}

Now we will consider more generally the moduli problem for families of marked complex analytic curves of genus $g$ with 
a given fixed subgroup of the automorphism group. 

\ms\noindent
Given a Riemann surface $C$ of genus $g$, let $G$ be a subgroup of its automorphism group $\Aut(C)$. We can view $G$ as a subgroup of the group of the 
outer automorphisms of $\Pi_g$, $\Out({\Pi}_g) = \Aut({\Pi}_g)/\Inn({\Pi}_g)$, since every complex automorphism which acts as the identity on the 
first homology group of $C$ must be the identity.

\ms\noindent
If we now look at a finite subgroup $G$ of the mapping class group $\Gamma_g = \pi_0(\Diff^+(\Sigma_g))$, as a subgroup of $\Out({\Pi}_g)$, we 
can consider the fixed point locus ${T}_g^G$ of the action of $G$ on the Teichm\"uller space ${T}_g$.

\ms\noindent
The complex analytic space ${T}_g^G$ is a nonempty complex submanifold of the Teichm\"uller space $T_g$ and it is also contractible. 
This follows from the solution of the Nielsen realization problem first proven by Kerkhoff \cite{K} and in the following form due to 
Tromba \cite{Tr} and slightly generalized by Catanese \cite{Ca1}.

\ms
\begin{proposition}
Given a finite subgroup $G$ of the mapping class group $\pi_0(\Diff(\Sigma_g))$, the fixed point locus ${T}_g^G$ of the action of $G$ on the 
Teichm\"uller space $T_g$ is nonempty and connected, indeed diffeomorphic to an euclidean space and therefore contractible.
\end{proposition}

\ms\noindent
We observe that ${T}_g^G$ parametrizes the isomorphism classes of Riemann surfaces of genus $g$ with a holomorphic action of $G$,
which is differentiably equivalent to the given one on the differentiable model $\Sigma_g$ of the Riemann surface $C$.

\ms\noindent
The submanifold $T_g^G$ is itself a Teichm\"uller space. In fact, let $\Sigma_g$ be a differentiable model for a fixed Riemann surface $C$ of 
genus $g$ and let us consider the quotient $\Sigma_g/G$ by the action of $G$. This quotient $\Sigma_{g'}:= \Sigma_g/G$ is a Riemann surface of 
a genus $g'$ and the projection $p: \Sigma_g \rightarrow \Sigma_{g'}$ is a ramified covering. We denote by $B = \{p_1,...,p_r\} \subset \Sigma_{g'}$ 
the branch locus of $p$. Let $T_{g',r}$ now be the Teichm\"uller space of the punctured surface $\Sigma^* := \Sigma_{g'} -B$.

\ms\noindent
A point in $T_g^G$ is just an isomorphism class $[(C, \phi, H)]$, where $\phi$ is the marking of the Riemann surface $C$ and $H=\nu(G)$ for the
given embedding $\nu: G \hookrightarrow \Aut(C)$ of $G$ in $\Aut(C)$.
For each such point $[(C, \phi, H)] \in T_g^G$, we denote by $C' = C/H$ the quotient under the action of $G$ and by $C^*$ the punctured surface 
obtained from $C'$ by removing the branch locus of the projection $C \rightarrow C'$. 
The marking $\phi: \pi_1(C) \stackrel{\cong} \rightarrow \Pi_g = \pi_1(\Sigma_g)$ induces a marking 
${\phi}^*: \pi_1(C^*) \rightarrow \pi_1(\Sigma^*)$, which defines a map 
$$\Psi: T_g^G \rightarrow T_{g',r}, \;\;\; [(C, \phi, H)] \mapsto [(C^*, {\phi}^*)].$$
We have the following fundamental theorem (cf. \cite{GDH}, \cite{H}, \cite{Kr}).

\ms
\begin{theorem}
\label{GDH}
The spaces $T_g^G$ and $T_{g',r}$ are biholomorphically equivalent via the mapping $\Psi$.
\end{theorem}

\ms\noindent
We also like to remind that any finite group $G$ actually has an embedding in the mapping class group of a certain Riemann surface.

\ms
\begin{proposition}
Given a finite group $G$, there is an embedding of $G$ in a mapping class group $\Gamma_g$ of a Riemann surface $C$ with genus $g$ 
for some $g \geq 2$. 
\end{proposition}

\noindent {\bf Proof.} Let us assume that the group $G$ is generated by $t$ elements $\delta_1,...,\delta_t$ and choose an integer $g' \geq 2$ such that $g' \geq t$. 
Let us define a homomorphism $\psi$ from the group $\Pi_{g'}$ to the free group $F_{g'}$ with $g'$ generators $\gamma_1, ..., \gamma_{g'}$ 
defined by $\gamma_i:=\psi(\alpha_i ) = \psi(\beta_i)$ for $i =1,...,g'$. 

\ms\noindent
Let us finally define an epimorphism $\varphi: F_{g'} \rightarrow G$ by sending $\gamma_i$ to $\delta_i$ for $i =1,...,t$, 
and such that $\varphi(\gamma_j)$ for $j >t$ can be chosen whatever we want. 

\ms\noindent
The composition $\varphi \circ \psi: \Pi_{g'} \rightarrow G$ is then an epimorphism and if we choose a point $[C', \phi]$ in the Teichm\"uller space 
$T_{g'}$, i.e., we choose a marking $\phi: \pi_1(C') \cong \Pi_{g'}$, we can consider the unramified Galois covering $C$ 
associated to the kernel of the epimorphism  $\varphi \circ \psi$. The genus of the Riemann surface $C$ is then by the Hurwitz formula given as 
$g = |G| (g'-1)+1$ and the group $G$ is contained in the group of automorphisms $\Aut(C)$ of $C$ and so we have an embedding of $G$
in the mapping class group $\Gamma_{g}$. 
\qed 

\ms\noindent
Let us now consider an interesting class of examples.

\ms\noindent {\bf Example.}
(cf. \cite{Ca1})
Consider the set $B \subset {\mathbb P}^1_{\C}$ consisting of the three real points $B :=\{-1,0,1\}$. Let us choose the following generators 
$\alpha, \ \beta, \ \gamma$ of the fundamental group $\pi_1({\mathbb P}^1_{\C}-B, 2)$:

\ms\noindent
$\alpha$ goes from $2$ to $-1 -\epsilon$ along the real line and passes through $+\infty$, then it makes a full turn counterclockwise along the 
circumference with centre $-1$ and radius $\epsilon$, then goes back to $2$ along the same way on the real line.    

\ms\noindent
$\gamma$ goes from $2$ to $1 + \epsilon$ along the real line, then it makes a full turn counterclockwise along the circumference with centre $+1$ 
and radius $\epsilon$, then goes back to $2$ along the same way on the real line.    

\ms\noindent
$\beta$ goes from $2$ to $+1 +\epsilon$ along the real line, it makes a half turn counterclockwise around the circumference with centre $+1$ and radius 
$\epsilon$, then it goes on the real line reaching $+\epsilon$, it makes a full turn counterclockwise around the circumference with centre $0$ and 
radius $\epsilon$, it goes back to $1-\epsilon$ along the real line, it makes a half turn clockwise around the circumference with centre $+1$ 
and radius $\epsilon$ and finally it proceeds along the real line returning to $2$. 

\ms\noindent
Now it is immediate to see that $\alpha$ and $\gamma$ are free generators of $\pi_1({\mathbb P}^1_{\C}-B, 2)$ and $\alpha \beta \gamma = 1$. Therefore we 
have fixed an isomorphism of $\pi_1({\mathbb P}^1_{\C}-B, 2)$ with the group $K := <\alpha, \ \beta, \ \gamma \ | \ \alpha \beta \gamma =1>$.

\ms\noindent
Let now $G$ be a finite group generated by two elements $a, \ b$. We have a surjection $\pi: K \rightarrow G$ given by $\pi(\alpha) = a$, $\pi(\beta) = b$, 
and we define $c:= \pi(\gamma)$. Then it follows that $abc =1$. 

\ms\noindent
The choice of the isomorphism $\pi_1({\mathbb P}^1_{\C}-B, 2) \cong K $ and of the epimorphism $\pi: K \rightarrow G$ determines a Galois covering 
$C \rightarrow {\mathbb P}^1_{\C}$ branched on $B$ with Galois group $G$.

\ms\noindent
Therefore $G$ is contained in $\Aut(C)$ and we have an embedding of $G$ in $\Gamma_g$ where $g$ is the genus of $C$. So we can consider the fixed 
point set $T_g^G$ of the action of $G$ on the Teichm\"uller space $T_g$ and by theorem \ref{GDH} we know that $T_g^G$ is isomorphic to $T_{0,3}$ 
which consists of only one point.  

\ms\noindent
Observe that if we set $m, \ n, \ p$ to be the periods of the respective elements $a, \ b, \ c$, then by the Hurwitz formula we have: 
$$2g-2 = |G|(1 - 1/m -1/n - 1/p).$$      

\ms\noindent
Now if we consider for instance the action on $C$ of the cyclic group $H$ generated by $a$, we find an embedding of $H$ in 
$\Gamma_g$ and consequently a fixed point locus $T_g^H$ of the action of $H$ on $T_g$.

\ms\noindent
An example of such a curve $C$ is the Fermat curve of degree $n$ in the projective space ${\mathbb P}^2_{\C}$. 

\ms\noindent
$C = \{(z_0,z_1,z_2) \ | \ z_0^n + z_1^n + z_2^n =0 \}$ is a Galois cover of the projective line ${\mathbb P}^1_{\C} = \{(x_0,x_1,x_2) \ | \ x_0 + x_1 + x_2 =0 \}$ under 
the map $\pi: {\mathbb P}^2_{\C}\rightarrow {\mathbb P}^2_{\C}$, given by $\pi(z_0,z_1,z_2) = (z_0^n,z_1^n, z_2^n)$. The map $\pi$ is ramified in three points 
and the Galois group of the covering is the group $G = ({\Z}/n)^2$ of diagonal projectivities with entries the $n$-th roots of unity. 

\ms\noindent
As it is observed in \cite{Ca1}, the curve $C$ can be seen as a covering of ${\mathbb P}^1_{\C}$ branched in three points in two different ways. 
In fact, we can consider the quotient of $C$ by the full group of automorphisms of $C$, which is a semidirect product of the normal subgroup $G$ 
by the symmetric group exchanging the three coordinates. 

\ms\noindent
This can be easily seen directly, but we also notice that if $G$ and $G'$ are two finite groups with $G \subset G' \subset  \Gamma_g$, we have an 
inclusion of the fixed point loci in the Teichm\"uller space $T_g$ as $T_g^{G'} \subset T_g^G$. 

\ms\noindent
Therefore in the case of the Fermat curve, we have $G \subset \Aut(C) \subset \Gamma_g$, where $g = (n-1)(n-2)/2$ is the genus of $C$ and thus 
$T_g^{\Aut(C)} \subset T_g^G \cong T_{0,3}$, which is only a point, i.e., the curve $C$ with the chosen marking, and since they are both nonempty we have 
$T_g^{\Aut(C)} = T_g^G \cong T_{0,3}$. 
So $C/\Aut(C)$ is isomorphic to ${\mathbb P}^1_{\C}$ and the map $C \rightarrow C/\Aut(C)$ is branched in three points.

\ms\noindent 
General actions of the group ${\Z}^m_p$ on Riemann surfaces were recently studied systematically also by Costa and Natanzon \cite{CN}.

\ms\noindent
Now we will introduce the moduli functor for families of marked complex analytic curves of genus $g\geq 2$ with a given 
fixed finite subgroup of the automorphism group as symmetry group.

\ms\noindent
Let $G$ be a finite group and let us assume that there exists an embedding of $G$ in the mapping class group $\Gamma_g$, i.e., there is an injection
$i:G \hookrightarrow \Gamma_g$. 

\ms
\begin{definition} 
A {\it family of marked complex analytic curves of genus $g$ with symmetry group $G$} is a triple $(p:C\rightarrow S, \phi, \nu: G \rightarrow \Aut_S(C))$,
where $p: C\rightarrow S$ is a family of complex analytic curves of genus $g$ over a complex analytic space $S$, $\phi$ a locally constant family of
markings of $p$ and $\nu: G \rightarrow \Aut_S(C)$ a group homomorphism injective in the fibres, with $G$ acting in such a way that for all points $s$
of $S$, if we denote by $\nu_s(G)_*$ the induced action of $\nu(G)$ on the fundamental 
group of the fibre $C_s$, and by $\psi: \Aut(\Pi_g) \rightarrow \Out(\Pi_g)$ the natural projection, we have  
$\psi(\phi_s \circ \nu_s(G)_* \circ \phi_s^{-1}) = i(G) \subset \Gamma_g$.
\end{definition}

\ms\noindent
To be more precise, for every element $h \in G$ and every point $s$ of $S$, $\nu_s(h)$ yields an automorphism of the fibre $C_s$, and therefore we have an 
isomorphism $\nu_s(h)_*: \pi_1(C_s, c_0) \stackrel{\cong} \rightarrow \pi_1(C_s, \nu_s(h)(c_0))$, for a basepoint $c_0$. 
Let us fix a path $\gamma$ between $\nu_s(h)(c_0)$ and $c_0$ and the corresponding isomorphism of fundamental groups $f_{\gamma}:\pi_1(C_s, \nu_s(h)(c_0)) 
\stackrel{\cong} \rightarrow \pi_1(C_s, c_0)$ 
sending an element $\delta$ to $\gamma^{-1} \delta \gamma$. 
Then if we take a representative $\phi_s$ of a marking $\phi$, say $\phi_s : \pi_1(C_s, c_0) \stackrel{\cong} \rightarrow \Pi_g$, 
we require in the definition that $\psi(\phi_s \circ f_{\gamma} \circ \nu_s(h)_* \circ \phi_s^{-1}) = i(h')$ for some element $h' \in G$. 
Clearly this does not depend on the chosen representative of the marking. 
Furthermore if we choose another path $\gamma'$ between $\nu_s(h)(c_0)$ and $c_0$ and the corresponding isomorphism of 
fundamental groups $f_{\gamma'}$, we see that 
$$(\phi_s \circ f_{\gamma'} \circ \nu_s(h)_* \circ \phi_s^{-1}) \circ (\phi_s \circ f_{\gamma} \circ \nu_s(h)_* \circ 
\phi_s^{-1})^{-1} = \phi_s \circ (f_{\gamma'} \circ f_{\gamma}^{-1}) \circ \phi_s^{-1}$$ 
is an element of $\Inn(\Pi_g)$, and so everything is well defined. 

\ms\noindent     
We are now going to define the moduli functor ${\mathcal T}_g^G$ from the category of analytic spaces to the category of sets as follows.

\ms
\begin{definition}
Let ${\mathcal T}_g^G: (AnSp)\rightarrow (Sets)$ be the contravariant functor from the category of complex analytic spaces to the
category of sets, which associates to every complex analytic space $S$ the set of isomorphism classes of families 
$(p:C\rightarrow S, \phi, \nu: G \rightarrow \Aut_S(C))$ of marked complex analytic curves with symmetry group $G$ over $S$, where isomorphism
between the triples $(p, \phi, \nu : G \rightarrow \Aut_S(C))$ and $(p', \phi', \nu' : G \rightarrow \Aut_S(C'))$ means isomorphism 
of the couples $(p, \phi)$ and $(p', \phi')$ (we do not require them to be equivariant). To every
holomorphic morphism $f: S'\rightarrow S$ the functor ${\mathcal T}_g^G$ associates the map between isomorphism classes 
induced by the base change with $f$.
\end{definition}

\ms\noindent
The fundamental statement is now, that the moduli problem for families of marked complex analytic curves of genus $g$ with symmetry
group $G$ has a fine solution, the fixed point locus of the Teichm\"uller space.

\ms\noindent
\begin{theorem}
The moduli functor ${\mathcal T}_g^G$ is representable by the complex analytic space $T_g^G$. 
\end{theorem}
\noindent {\bf Proof.} The fixed point locus $T_g^G$ is a complex submanifold of $T_g$ and we claim that if $\alpha
: {\mathcal U}_g \rightarrow T_g$ is the universal family on $T_g$ and $j:
T_g^G \hookrightarrow T_g$ the embedding of $T_g^G$ in $T_g$, then the
fibre product $T_g^G \times_{T_g} {\cal U}_g =:{\cal U}_g^G$ defines a universal family on $T_g^G$.

\ms\noindent
Since $j: T_g^G \hookrightarrow T_g$ is an embedding, then also
the induced map $f: {\cal U}_g^G \hookrightarrow {\cal U}_g$ is an embedding. 

\ms\noindent
Now, given a family $p: C \rightarrow S$ together with a marking $\phi$ and a homomorphism $\nu: G
\rightarrow \Aut_S({\mathcal C})$ representing an isomorphism class $(p, \phi, \nu)$ in ${\mathcal T}_g^G(S)$, 
we have a map $\gamma : S \rightarrow T_g$ obtained by
forgetting the action of the group $G$. As $\gamma$ factors through 
$T_g^G$ via a map $h: S \rightarrow T_g^G$ we have $\gamma = j\circ h$. 

\ms\noindent
Then, since $\alpha: {\cal U}_g \rightarrow T_g$ is a universal
family, we know that $C$ is given as the fibre product $S \times_{T_g} {\cal U}_g$, and we get a cartesian diagram

$$\diagram
   C \rto^{\beta} \dto_{p}&               {\cal U}_g\dto^{\alpha} \\
   S  \rto^{\gamma}&        T_g               \enddiagram $$

\noindent where $\beta$ factors through a map $k : C \rightarrow {\cal U}_g^G$. Also the following diagram 

$$\diagram
   {\cal U}_g^G \rto^{f} \dto &               {\cal U}_g\dto^{\alpha} \\
   T_g^G  \rto^{j}&        T_g               \enddiagram $$ 

\noindent is cartesian and so we get a big commutative diagram 

$$\diagram
  C\rto^{k}\dto_{p} & {\cal U}_g^G \rto^{f} \dto &               {\cal U}_g\dto^{\alpha} \\
  S\rto^{h}&           T_g^G  \rto^{j}&        T_g               \enddiagram $$ 

\noindent such that the diagram on the right is cartesian and the one obtained
by taking the composition maps is also cartesian. Since the map $f : {\cal
U}_g^G \hookrightarrow {\cal U}_g$ is injective, it is immediate to
verify the universal property in order to show that also the diagram on the left

$$\diagram
  C\rto^{k}\dto_{p} & {\cal U}_g^G \dto  \\
  S\rto^{h} &           T_g^G     \enddiagram $$ 

\noindent is cartesian.

\ms\noindent
The uniqueness of the map $h$ is a direct consequence of the universality of the family ${\cal U}_g \rightarrow T_g$. 

\ms\noindent
So we have proven that ${\cal U}_g^G \rightarrow T_g^G$
is a universal family and therefore the moduli functor ${\mathcal T}_g^G$ is representable by the complex analytic space $T_g^G$.
\qed

\ms\noindent
Because of the representability of the moduli functor ${\mathcal T}_g^G$, the identity morphism $id: T_g^G\rightarrow T_g^G$ 
defines the universal family of marked complex analytic curves with symmetry group $G$ over $T_g^G$, which, as the
proof above shows, is given as a triple $({\cal U}_g^G\rightarrow T_g^G, \Phi, \nu: G\rightarrow \Aut_{T_g^G}({\cal U}_g^G))$,
where $T_g^G$ is the fixed point set of the classical Teichm\"uller space $\T_g$ and ${\cal U}_g^G$ the fibre product $T_g^G\times_{T_g} {\cal U}_g$
along the classical universal family ${\cal U}_g\rightarrow T_g$ given by the identity morphism $id: T_g\rightarrow T_g$ via
the representability of the moduli functor ${\cal T}_g$ (cf. Theorem 2.5).

\subsection{The moduli stack of complex analytic curves with symmetries and its analytic homotopy type}

\ms
We will now study the moduli stack of families of complex analytic curves of genus $g$ with a fixed finite group $G$ of automorphisms
and determine its analytic homotopy type as introduced in the first chapter. We will always assume here again that $g\geq 2$. Forgetting
markings we get the concept of families of complex analytic curves with symmetries.

\ms
\begin{definition}
Let $S$ be a complex analytic space and $G$ a finite group with fixed embedding in the mapping class group $\Gamma_g$.
A {\it family of complex analytic curves of genus $g$ over $S$ with symmetry group $G$} is a pair $(p:C\rightarrow S, \nu: G\rightarrow \Aut_S(C))$,
where $p: C\rightarrow S$ is a family of complex analytic curves of genus $g$ over $S$ and $\nu: G\rightarrow \Aut_S(C)$ is a group homomorphism
such that the group $\nu(G)$ of $S$-automorphisms is acting faithfully on the fibres $C_s$ and this action is differentiably equivalent to the
one given by the embedding of $G$ in $\Gamma_g$.
\end{definition}

\ms\noindent
We define now the analytic moduli stack ${\mathcal M}_{g,G}^{an}$ of families of complex analytic curves of genus $g$ with symmetries.

\ms
\begin{definition}
Let $G$ be a finite group embedded in the mapping class group $\Gamma_g$. 
The moduli stack ${\mathcal M}_{g,G}^{an}$ of {\it families of complex analytic curves of genus $g$ with symmetrie group G} is the category fibred in groupoids 
over the category (AnSp) of complex analytic spaces defined by its groupoid of sections ${\mathcal M}_{g,G}^{an}(S)$ as follows: For a complex analytic space $S$ 
the objects of ${\mathcal M}_{g,G}^{an}(S)$ are families $(p:C\rightarrow S, \nu: G\rightarrow \Aut_S(C))$ of complex analytic curves of genus $g$ over $S$ with 
symmetry group $G$ and the morphisms are $G$-equivariant $S$-isomorphisms, i.e., morphisms between $(p: C \rightarrow S, \ \nu:G \rightarrow \Aut_S(C))$ and
$(p': C' \rightarrow S, \  \nu':G \rightarrow  \Aut_S(C'))$ are $S$-isomorphisms $\phi: C \rightarrow C'$ 
such that $\nu(G) = \phi^{-1} \circ \nu'(G) \circ \phi$.
\end{definition}

\ms\noindent
Consider now the universal family $({\cal U}_g^G\rightarrow T_g^G, \Phi, \nu: G\rightarrow \Aut_{T_g^G}({\cal U}_g^G))$ on $T_g^G$, where $\Phi$ 
is the universal marking, then by forgetting the universal marking, this canonical data defines an object in the category of 
sections ${\mathcal M}_{g,G}^{an}(T_g^G)$ of the stack ${\mathcal M}_{g,G}^{an}$.

\ms\noindent
Considering the complex analytic space $T_g^G$ as a stack, i.e., as the stack associated to the complex analytic space $T_g^G$, 
and forgetting the universal marking induces a morphism of stacks $\pi: T_g^G \rightarrow {\mathcal M}_{g,G}^{an}$.
It is not hard to see from this that ${\mathcal M}_{g,G}^{an}$ is actually an analytic Deligne-Mumford stack as
introduced in the first chapter. As we shall see in the next chapter ${\mathcal M}_{g,G}^{an}$ is in fact the complex 
analytification of an algebraic Deligne-Mumford stack.

\ms\noindent
Now we determine the fibre product $Isom(\pi, \pi):=T_g^G \times_{{\mathcal M}_{g,G}^{an}} T_g^G$ of the morphism $\pi$ of stacks with itself.
We have the following diagram

$$\diagram
Isom(\pi, \pi)= T_g^G \times_{{\mathcal M}_{g,G}^{an}} T_g^G \rto \dto &     T_g^G \dto^{\pi} \\
  T_g^G \rto^{\pi} &            {\mathcal M}_{g,G}^{an}  \enddiagram $$ 

\noindent Both projections from $Isom(\pi, \pi)$ to $T_g^G$ are local isomorphisms for the classical topology (cf. \cite{M}).

\ms\noindent
For each point $x \in T_g^G$ the data over $x$ has no nontrivial automorphism, therefore the map 
$Isom(\pi, \pi) \rightarrow T_g^G \times T_g^G$ is actually an immersion.

\ms\noindent
We know that two points $x,y \in T_g$ define the same object over the stack ${\mathcal M}_g^{an}$ if and only if the fibres $C_x$ and $C_y$ are 
isomorphic curves and for the markings $\Phi_x$ of $C_x$ and $\Phi_y$ of $C_y$ we have that $\Phi_y = \theta \circ \Phi_x$, where $[\theta]:=\psi(\theta)$ is an element
of the mapping class group $\Gamma_g = A(\Pi_g)/\Inn(\Pi_g)$. (cf. \cite{O}). Here $A(\Pi_g)$ is the subgroup of $\Aut(\Pi_g)$ given by the automorphisms 
inducing the identity in the second cohomology $H^2(\Pi_g, \Z) \cong \Z$. 

\ms\noindent
If now two points $x, y \in T_g^G$ i.e., two isomorphism classes of marked complex analytic curves of genus $g$ with symmetry group $G$, define 
the same object over the stack ${\mathcal M}_{g,G}^{an}$, there must exists an isomorphism of the curves 
$h: C_x \stackrel{\cong} \rightarrow  C_y$ which conjugates $\nu_x(G)$ and $\nu_y(G)$, i.e., $h^{-1}  \nu_y(G) h = \nu_x(G)$, where
$\nu_x: G\rightarrow \Aut(C_x)$ and $\nu_y: G\rightarrow \Aut(C_y)$ are the induced actions of $G$ on the fibres $C_x$ and $C_y$ and furthermore 
$\psi(\Phi_x (\nu_x (G))_* \Phi_x^{-1}) = i(G)$ and $\psi(\Phi_y (\nu_y (G))_*   \Phi_y^{-1}) = i(G)$, where $i:G\rightarrow \Gamma_g$ is the
fixed embedding of the group $G$ in the mapping class group $\Gamma_g$.

\ms\noindent
Since $\theta = \Phi_y \circ h_* \circ \Phi_x^{-1}$, we get immediately
$$\theta \psi^{-1}(i(G)) \theta^{-1} = \Phi_y h_* \Phi_x^{-1} \psi^{-1}(i(G)) \Phi_x h_*^{-1} \Phi_y^{-1} = \Phi_y h_* (\nu_x(G))_* h_*^{-1} \Phi_y^{-1} = $$
$$= \Phi_y (\nu_y(G))_* \Phi_y^{-1} = \psi^{-1}(i(G)).$$
And so we finally have
$$[\theta] i(G) [\theta]^{-1}= i(G).$$ 

\ms\noindent
Therefore we have proven that for each point $x \in T_g^G$, and if we identify $G$ with $i(G) \subset \Gamma_g$, the fibre $p^{-1}(x)$ of the 
projection $p: Isom(\pi, \pi) \rightarrow T_g^G$ is isomorphic to the cartesian product $T_g^G \times \Gamma_{g,G}$, where the
group $\Gamma_{g, G} := \{\alpha \in \Gamma_g \ | \ \alpha G \alpha^{-1} = G, \}$ is the normalizer of the finite group $G$ in the mapping class group $\Gamma_g$.
In fact, the isomorphism is explicitely given as
$$p^{-1}(x) \stackrel{\cong} \rightarrow T_g^G \times \Gamma_{g, G},$$
$$(x,y) \mapsto (x, [\Phi_y \circ \Phi_x^{-1}] = [\theta]).$$
By induction on the number of factors, we easily get an isomorphism
$$T_g^G \times_{{\mathcal M}_{g,G}^{an}}T_g^G \hdots \times_{{\mathcal M}_{g,G}^{an}} T_g^G \cong T_g^G \times \Gamma_{g, G} \times \hdots \times \Gamma_{g, G}.$$

\ms\noindent
Now we determine the analytic homotopy type of the moduli stack ${\mathcal M}_{g,G}^{an}$ of complex analytic curves with 
symmetries.

\ms
\begin{proposition}
Let $\cosk_0^{{\mathcal M}_{g,G}^{an}}(T_g^G)$ be the \v{C}ech nerve associated to the locally
isomorphic surjective covering morphism $\pi: T_g^G \rightarrow {\mathcal M}_{g,G}^{an}$.
Then its geometric realization $|\cosk_0^{{\mathcal M}_{g,G}^{an}}(T_g^G)|$ is homotopy
equivalent to the classifying space $|B\Gamma_{g,G}|$ of the group $\Gamma_{g, G}$.
\end{proposition}

\ms\noindent
{\bf Proof.} From the induction argument above, we see first that the $m$-simplex $\cosk_0^{{\mathcal M}_{g,G}^{an}}(T_g^G)_m$ of the
\v{C}ech nerve $\cosk_0^{{\mathcal M}_{g,G}^{an}}(T_g^G)$ is given by the $(m +1)$-tuple
fiber product of copies of $T_g^G$ over ${\mathcal M}_{g,G}^{an}$, which as we showed is isomorphic
to $T_g^G \times \Gamma_{g,G} \times \Gamma_{g, G} \times \cdots \Gamma_{g, G}$. So by definition, we get
$$\cosk_0^{{\mathcal M}_{g,G}^{an}}(T_g^G)\cong T_g^G \times B\Gamma_{g,G},$$
where $B\Gamma_{g, G}$ is the classifying simplicial set of the discrete
group $\Gamma_{g, G}$. Therefore after geometric realization we have a homeomorphism of topological spaces
$$|\cosk_0^{{\mathcal M}_{g,G}^{an}}(T_g^G)|\cong T_g^G \times |B\Gamma_{g, G}|.$$
From Proposition 2.6 we know that the space $T_g^G$ as an euclidean space
is indeed contractible, so we finally get the conclusion of the proposition.
\qed

\ms\noindent
The following corollary is therefore an immediate consequence from the considerations
about analytic homotopy types of the first chapter as the \v{C}ech nerve $\cosk_0^{{\mathcal M}_{g,G}^{an}}(T_g^G)$
is a simplicial analytic space and a hypercovering of the topos $\hf\vf ({{\mathcal M}_{g,G}^{an}}_{cl})$.

\ms
\begin{corollary}
There is a weak homotopy equivalence of simplicial sets
$$ \{{\mathcal M}_{g,G}^{an}\}_{cl}\simeq \Sin(|B\Gamma_{g, G}|).$$

\end{corollary}

\ms\noindent
Let us introduce another analytic moduli stack ${\mathcal M'}_{g,G}^{an}$, which is closely related to
the moduli stack ${\mathcal M}_{g,G}^{an}$ (cf. also \cite{GDH}).

\ms
\begin{definition}
Let $G$ be a finite group embedded in the mapping class group $\Gamma_g$. 
The moduli stack ${\mathcal M'}_{g,G}^{an}$ is the category fibred in groupoids over the category $(AnSp)$ of complex analytic spaces defined
by its groupoid of sections ${\mathcal M'}_{g,G}^{an}(S)$ as follows: For a complex analytic space $S$ the objects of ${\mathcal M'}_{g,G}^{an}(S)$
are families $(p:C\rightarrow S, \nu: G\rightarrow \Aut_S(C))$ of complex analytic curves of genus $g$ over $S$ with symmetry group $G$ and 
the morphisms are $S$-isomorphisms, not necessarily equivariant, i.e., morphisms between $(p: C \rightarrow S, \ \nu:G \rightarrow \Aut_S(C))$ and
$(p': C' \rightarrow S, \  \nu':G \rightarrow  \Aut_S(C'))$ are just $S$-isomorphisms $\phi: C \rightarrow C'$. 
\end{definition}

\ms\noindent
While the moduli stack ${\mathcal M}_{g,G}^{an}(S)$ is not an analytic substack of ${\mathcal M}^{an}_g$, the moduli stack ${\mathcal M'}_{g,G}^{an}$ 
is and by forgetting the action we get a natural morphism of stacks $f: {\mathcal M}_{g,G}^{an} \rightarrow {\mathcal M'}_{g,G}^{an}$.
We also have the following commutative diagram:

$$\diagram
T_g^G \dto^{\pi} \rrto &   &  T_g \dto \\
  {\mathcal M}_{g,G}^{an} \rto^{f} &  {\mathcal M'}_{g,G}^{an} \rto &  {\mathcal M}_{g}^{an} \enddiagram $$ 

\noindent where the map $\pi$ is the morphism of stacks as defined above.

\ms\noindent
Similarly as above we can see that two points $x,y \in T_g^G$, i.e., two isomorphism classes of marked complex analytic curves of genus $g$ with
symmetry group $G$
$(C_x, \Phi_x, \nu_x: G \hookrightarrow \Aut(C_x))$ and $(C_y, \Phi_y, \nu_y: G \hookrightarrow \Aut(C_y))$ give the same image in the stack ${\mathcal M'}^{an}_{g,G}$ 
if and only if $C_x \cong C_y$ and $\Phi_y = \theta \circ \Phi_x$ where $[\theta] \in \Gamma_g$. 
But we also must have that $(\Phi_x)^{-1} \circ G \circ \Phi_x = (\nu_x (G))_* =:K$ and  
$(\Phi_y)^{-1} \circ G \circ \Phi_y = (\nu_y (G))_*  = (\theta \circ \Phi_x)^{-1} \circ G \circ (\theta \circ \Phi_x) =: K'$, where $K, K'$ are
subgroups of $\Aut(C_x)$, since we have identified $C_x$ with $C_y$.  Now if $(C_x, \Phi_x,\nu_x)$ and $(C_y, \Phi_y,\nu_y)$ give different images 
in ${\mathcal M}^{an}_{g,G}$, then $[\theta] \not\in \Gamma_{g,G}$ and so $[\theta]^{-1} G [\theta] \neq G$, hence $K \neq K'$ are two different 
subgroups of $\Aut(C_x)$.

\ms\noindent 
This shows therefore that the moduli stacks ${\mathcal M}^{an}_{g,G}$ and ${\mathcal M'}^{an}_{g,G}$ are different if and only if there 
exists a Riemann surface $C$ of genus $g$ with symmetry group $G$ whose automorphism group contains two subgroups which are 
conjugated topologically, but not holomorphically. 

\ms\noindent From the results of Gonz\'alez-D\'{\i}ez and Harvey \cite{GDH} it also follows that ${\mathcal M}^{an}_{g,G}$ can actually be described
as an orbifold given as the quotient $T_g^G/\Gamma_{g,G}$, where $T_g^G$ is the fixed point locus of the classical Teichm\"uller space $T_g$ 
and where the modular group $\Gamma_{g,G}$ is the normalizer of the embedded finite group $G$ in the mapping class group $\Gamma_g$. Actually 
${\mathcal M}^{an}_{g,G}$ is a normal complex analytic space, the normalization of the irreducible subvariety ${\mathcal M'}^{an}_{g,G}$ of ${\mathcal M}^{an}_g$.

\ms\noindent
\begin{example}
If we take the group $G = {\Z}/2$ and $C$ is a hyperelliptic complex analytic curve with hyperelliptic involution
$\tau: C \rightarrow C$, then $\tau$ is the only automorphism of the curve $C$
of order two such that the quotient $C/<\tau>$ is isomorphic to the complex projective line ${\mathbb{P}}^1(\C)$, therefore
we have ${\mathcal M}^{an}_{g,<\tau>} = {\mathcal M'}^{an}_{g,<\tau>}$ and the moduli stack ${\mathcal M}^{an}_{g,<\tau>}$ of hyperelliptic complex analytic
curves of genus $g$ is actually an analytic substack of the moduli stack ${\mathcal M}_g^{an}$.
\end{example}

\section{Moduli stacks of families of algebraic curves with symmetries and their etale homotopy types}

\subsection{The Deligne-Mumford stack of algebraic curves with symmetries}

\ms
In this section we will study moduli stacks of algebraic curves with symmetries. Let $(Sch)$ here always be the category of
schemes over $\Q$. First we recall some basic notions (cf. \cite{MFK}, \cite{DM}).

\ms
\begin{definition}
Let $g\geq 2$ and $S$ be a scheme (over $\Q$). A {\it family of algebraic curves of genus $g$ over $S$} or {\it algebraic curve over $S$}
is a morphism $p:C\rightarrow S$ of schemes such that $p$ is proper and smooth and the geometric fibers $C_s$ of $p$ are 
$1$-dimensional smooth connected schemes of genus $g$. We will also denote an algebraic curve over $S$ simply by $C/S$.
\end{definition}

\ms\noindent
Here geometric fibers means as usual scheme theoretic fibers over points, i.e
for a point $s$ of $S$ we define $C_s:=C\times_S \Spec \,k(s)$, with $k(s)=\O_{S,s}/{\mathfrak{m}}_{S,s}$. 
And the genus is given cohomologically as $g=\dim H^1(C_s, \O_{C_s})$.
Over the complex numbers $\C$ geometric fibers are just complex smooth projective curves and morphisms are
regular maps and in the associated complex analytic case geometric fibers are just compact Riemann surfaces of topological 
genus $g$.

\ms
\begin{definition}
Let $S$ be a scheme (over $\Q$) and $G$ a finite group. An {\it action of the group $G$} on a family of algebraic curves $p: C \rightarrow S$ over $S$ 
is a morphism of group schemes over $S$ 
$$\mu: G_S\rightarrow \Aut_S(C)$$
of the constant group scheme $G_S: = G$ over $S$ to $\Aut_S(C)$, which is injective in the fibres $C_s$.
\end{definition}

\ms\noindent
Let us now again consider a finite group $G$ with a fixed embedding $i:G \hookrightarrow \Gamma_g$ in the mapping class group $\Gamma_g$. 
We define now the moduli stack ${\mathcal M}_{g,G}$ of algebraic curves with symmetries.

\ms
\begin{definition}
The moduli stack ${\mathcal M}_{g,G}$ of {\it algebraic curves with symmetries} is the category fibred in groupoids over the category $(Sch)$ of schemes over $\Q$ defined
by its groupoid of sections ${\mathcal M}_{g,G}(S)$ as follows:
For a scheme $S$ over $\Q$ the objects of ${\mathcal M}_{g,G}(S)$ are families of algebraic curves $p: C \rightarrow S$ over $S$ such that the fibres 
are connected smooth algebraic curves of genus $g \geq 2$ endowed with an action of the group $G$ on $p: C \rightarrow S$ satisfying the following 
property: Let $C_s$ be the fibre over a point $s$ of $S$. 
It is a scheme over $\Spec(\Q)$ and we let $C_s^{an}$ to be the complex analytic space associated with the $\C$-valued points $C_s(\C)$ of $C_s$. 
We require that the differentiable action of $G$ on $C_s^{an}$ is the one given by the embedding of $G$ in $\Gamma_g$ that we have fixed.
The morphisms are $G$-equivariant $S$-isomorphisms, i.e., morphisms between $(p: C \rightarrow S, \nu:G \rightarrow \Aut_S(C))$ and 
$(p': C' \rightarrow S, \  \nu':G \rightarrow  \Aut_S(C'))$
are $S$-isomorphisms $\phi: C \rightarrow C'$ such that $\nu(G) = \phi^{-1} \circ \nu'(G) \circ \phi$.
\end{definition}

\ms\noindent
We will also consider the following moduli stack, studied by Tuff\'ery \cite{Tu1}, \cite{Tu} which we will
compare with ${\mathcal M}_{g,G}$ via forgetting the embedding of the finite group $G$ in the mapping class group $\Gamma_g$.  

\ms
\begin{definition}
Let $G$ be a finite group. The moduli stack ${\mathcal M}_g[G]$ is the category fibred in groupoids over the category $(Sch)$ of schemes over $\Q$ defined
by its groupoid of sections ${\mathcal M}_g[G](S)$ as follows:
For a scheme $S$ over $\Q$ the objects of ${\mathcal M}_g[G](S)$ are families of algebraic curves $p: C \rightarrow S$ over $S$ 
with a group action of the group $G$ and the morphisms are the $G$-equivariant $S$-isomorphisms.
\end{definition}

\ms\noindent
The moduli stack ${\mathcal M}_{g,G}$ is a substack of ${\mathcal M}_g[G]$ and its main property is the following.

\ms
\begin{proposition}
The moduli stack ${\mathcal M}_{g,G}$ is a Deligne-Mumford stack.
\end{proposition}

\ms\noindent
{\bf Proof.} Let ${\mathcal M}_{g} \otimes \Q$ be the moduli stack which represents the restriction of the Deligne-Mumford stack ${\mathcal M}_{g}$ of
all families of algebraic curves of genus $g\geq 2$ (cf. \cite{DM}) to the subcategory of schemes over $\Q$.

\ms\noindent  
We claim that the morphism of stacks ${\mathcal M}_g[G] \rightarrow {\mathcal M}_{g} \otimes \Q$ induced by forgetting the action of the group 
$G$ is finite and representable (cf. \cite{Tu}), thus the stack ${\mathcal M}_g[G]$ is a Deligne-Mumford stack and the same holds then for 
${\mathcal M}_{g,G}$ as being a substack of ${\mathcal M}_{g}[G]$. 

\ms\noindent
Recall that a morphism ${\mathcal{N}}\rightarrow \mathcal{M}$ of stacks is said to be representable if for every morphism $Y \rightarrow \mathcal{M}$ of
stacks,  where $Y$ is (the stack associated to) a scheme, the fibre product $Y \times_{\mathcal{M}} \mathcal{N}$ is isomorphic to 
(the stack associated to) a scheme. 

\ms\noindent
Therefore, in order to prove that the morphism of stacks $g: {\mathcal M}_g[G] \rightarrow {\mathcal M}_{g} \otimes \Q$ given by forgetting the action 
of $G$ is representable, we must prove that given a morphism $f: Y \rightarrow {\mathcal M}_{g} \otimes \Q$, with $Y$ a scheme, any object in 
the fibre product $Z = Y \times_{{\mathcal M}_{g} \otimes \Q} {\mathcal M}_g[G]$ has no nontrivial automorphisms. 

\ms\noindent
The objects of $Z$ are triples $(x,z, \alpha)$ where $\alpha:f(x) \rightarrow g(z)$ is a morphism in a fibre of ${\mathcal M}_{g} \otimes \Q$. 
A morphism $(x,z, \alpha) \rightarrow (x',z', \alpha')$ is a pair of morphisms $(\beta_1: x \rightarrow x'; \beta_2: z \rightarrow z')$ in the 
fibres of $Y$ and ${\mathcal M}_g[G]$ respectively such that $g(\beta_2) \circ \alpha = \alpha' \circ f(\beta_1): f(x) \rightarrow g(z')$.

\ms\noindent
So, assume that an object $(x,z, \alpha)$ of $Z$ has an automorphism $(\beta_1: x \rightarrow x; \beta_2: z \rightarrow z)$ as above, 
then clearly $\beta_1 = id$, and we must have $g(\beta_2) \circ \alpha = \alpha $, so $g(\beta_2) = id$. 
Now, since $z$ is a family $C \rightarrow S$ of curves with a $G$-action and $\beta_2$ is an equivariant automorphism of 
$C \rightarrow S$, we have that $g(\beta_2) =\beta_2$, because $g$ only forgets the action of $G$. 
Thus we have shown $\beta_2 = g(\beta_2) = id$ and $(x,z, \alpha)$ has no nontrivial automorphism.

\ms\noindent
This proves that the stack $Z$ is indeed a scheme and the morphism of stacks $g: {\mathcal M}_g[G] \rightarrow {\mathcal M}_{g} \otimes \Q$ 
is a representable morphism. 
\qed

\subsection{The etale homotopy type of the moduli stack of algebraic curves with symmetries}

\ms
Now we can prove our main theorem on the etale homotopy type of the moduli stack of families of algebraic curves
of genus $g$ with a prescribed symmetry group using the above considerations and the general comparison theorem
for etale homotopy types of Cox and Friedlander as stated in the first section.

\ms\noindent
Let $G$ be again a finite group with an embedding in the mapping class group $\Gamma_g$ of genus $g$.
Let now ${\mathcal M}_{g,G} \otimes \bar{\Q}$ be the restriction of the moduli stack ${\mathcal M}_{g, G}$ to the subcategory of 
schemes over $\bar{\Q}$. Fix an embedding $\bar{\Q} \hookrightarrow {\C}$ of the algebraic closure of the rationals in the complex numbers.
The analytic stack ${\mathcal M}_{g,G}^{an}$ is precisely the complex analytification 
$({\mathcal M}_{g,G}\otimes\bar{\Q})^{an}$ of the Deligne-Mumford stack ${\mathcal M}_{g,G}\otimes\bar{\Q}$ and
we can now determine the etale homotopy type of ${\mathcal M}_{g,G}\otimes\bar{\Q}$.

\ms
\begin{theorem}
There is a weak homotopy equivalence of pro-simplicial sets
$$\{{\mathcal M}_{g, G}\otimes\bar{\Q}\}^{\wedge}_{et}\simeq
 K(\Gamma_{g,G}, 1)^{\wedge}.$$
where $\Gamma_{g, G}=\Norm_{\Gamma_g}(G)$ is the normalizer of the group $G$ in the mapping class group $\Gamma_g$.
\end{theorem}

\ms\noindent
{\bf Proof.} The moduli stack ${\mathcal M}_{g,G}\otimes\bar{\Q}$ is a Deligne-Mumford stack, so there exists a 
scheme $X$ of finite type over $\bar{\Q}$ and an etale surjective covering morphism
$x: X\rightarrow {\mathcal M}_{g,G}\otimes\bar{\Q}.$
The \v{C}ech nerve $\cosk_0^{{\mathcal M}_{g,G}\otimes\bar{\Q}}(X)$ for
this morphism defines a hypercovering of the stack ${\mathcal M}_{g,G}\otimes\bar{\Q}$.
Similarly $\cosk_0^{({\mathcal M}_{g,G}\otimes\bar{\Q})^{an}}(X^{an})$ is a hypercovering of the analytic
stack ${\mathcal M}_{g,G}^{an}$. Here $X^{an}$ denotes the associated complex analytic
space of the covering scheme $X$ over $\bar{\Q}$.

\ms\noindent
The homotopy descent theorem (Theorem 1.3) shows that there are weak equivalences of pro-simplicial sets:
$$\{\m_{g,G}\otimes\bar{\Q}\e \simeq \{\cosk_0^{\m_{g,G}\otimes\bar{\Q}}(X)\e \simeq 
\{\cosk_0^{\m_{g,G}\otimes\bar{\Q}}(X)\otimes_{\bar{\Q}}\C\e.$$

\ms\noindent
and the comparison theorem for simplicial schemes (cf. Cox \cite{C}, Theorem IV.8) shows that there is also a weak
equivalence after profinite completions:

$$\{\cosk_0^{\m_{g,G}\otimes\bar{\Q}}(X)\otimes_{\bar{\Q}}\C\}^{\wedge}_{et} \simeq 
\Sin(|\cosk_0^{({\mathcal M}_{g,G}\otimes\bar{\Q})^{an}}(X^{an})|)^\wedge.$$

\ms\noindent
As $\cosk_0^{({\mathcal M}_{g,G}\otimes\bar{\Q})^{an}}(X^{an})$ is a hypercovering of the analytic
stack ${\mathcal M}_{g,G}^{an}$ we have a weak equivalence:

$$\{\m_{g,G}^{an}\}_{cl} \simeq \Sin(|\cosk_0^{({\mathcal M}_{g,G}\otimes\bar{\Q})^{an}}(X^{an})|)$$

\ms\noindent
From the determination (Corollary 2.15) of the classical homotopy type of ${\mathcal M}_{g,G}^{an}$ we get therefore
at the end the following weak homotopy equivalence:
$$\{{\mathcal M}_{g, G}\otimes\bar{\Q}\}^{\wedge}_{et}\simeq
 B\Gamma^{\wedge}_{g, G} =K(\Gamma_{g,G}, 1)^{\wedge}.$$

\noindent This finally proves the main theorem. \qed

\ms\noindent
If the group $G$ is the trivial group we recover as a special case the main theorem of Oda (cf. \cite{O}, Theorem 1) for 
the moduli stack $\m_g$ of all families of algebraic curves of genus $g \geq 2$.

\ms\noindent
Another interesting special case is given by the moduli stack $\h_g$ parametrizing only families of hyperelliptic curves
of genus $g$ (cf. also \cite{LL} and \cite{LK}).

\ms
\begin{definition}
Let $g\geq 2$ and $S$ be a scheme (over $\Q$). A family of algebraic curves $p: C\rightarrow S$ of genus $g$  over $S$ is called {\it  hyperelliptic} if there
exists an $S$-involution $\tau: C\rightarrow C$ such that $C/<\tau>\cong {\mathbb{P}}^1_S$.
\end{definition}

\ms\noindent
We define the moduli stack ${\mathcal H}_g$ of hyperelliptic curves of genus $g$ as follows.

\ms
\begin{definition}
The moduli stack ${\mathcal H}_g$ of {\it hyperelliptic curves} is the category fibred in groupoids over the category $(Sch)$ of schemes over $\Q$ defined
by its groupoid of sections ${\mathcal H}_g(S)$ as follows:
For a scheme $S$ over $\Q$ the objects of ${\mathcal H}_g(S)$ are families of hyperelliptic curves $p: C \rightarrow S$ of genus $g$ 
over $S$ and the morphisms are $S$-isomorphisms equivariant with respect to the hyperelliptic involution.
\end{definition}


\ms\noindent
The moduli stack of hyperelliptic curves $\h_g$ is a Deligne-Mumford stack. In fact, $\h_g$ is also one of the {\it Hurwitz stacks} studied by Fulton \cite{Fu}, \S 6, namely $\h_g$ is the stack $\h^{2,2g +2}$ of 2-sheeted simple coverings of ${\mathbb{P}}^1_S$ with $2g +2$ 
branching points, where $S$ is a scheme (over $\Q$). It would be interesting to study the etale homotopy types of more general Hurwitz stacks.

\ms\noindent
In our words the stack $\h_g$ of hyperelliptic curves is simply the moduli stack $\m_{g,G}$, where $G=<\tau>$
is the group generated by a hyperelliptic involution $\tau$ as discussed above, i.e., $G\cong \Z/2$.

\ms\noindent
In this particular case the group $\Gamma_{g, G}$ is then just the {\it hyperelliptic mapping class group group} $\Gamma^h_g$,
which is the centralizer $\Cent_{\Gamma_g}(<\tau>)$ of a hyperelliptic involution $\tau$ in the mapping class 
group $\Gamma_g$.

\ms\noindent
It also turns out, as already mentioned by Oda \cite{O} that in this case the etale homotopy type 
of the moduli stack of hyperelliptic curves is an honest Eilenberg-MacLane space  $K(\Gamma^{h \wedge}_g, 1)$, where $\Gamma^{h \wedge}_g$ is the 
profinite completion of $\Gamma^h_g$. This is because the hyperelliptic mapping
class group $\Gamma^h_g$ is a {\em good} group in the sense of Serre \cite{S}, which essentially means that the Galois cohomology of its profinite 
completion coincides with the corresponding discrete group cohomology. Therefore there is a weak homotopy equivalence after profinite Artin-Mazur 
completion between $K(\Gamma^h_g, 1)$ and $K(\Gamma^{h \wedge}_g, 1)$ as was shown by Artin and Mazur \cite{AM}, 6.6. And we can deduce now directly 
from our main theorem the following corollary:

\ms
\begin{corollary}
Let $\h_g$ be the moduli stack of hyperelliptic curves of genus $g$. There is a weak homotopy equivalence of pro-simplicial sets
$$\{\h_g\otimes\bar{\Q}\}^{\wedge}_{et} \simeq K(\Gamma^{h \wedge}_g, 1).$$
where $\Gamma^h_g$ is the hyperelliptic mapping class group. 
\end{corollary}

\ms\noindent
From the discussion above, the following open question seems to be of some general interest:

\ms\noindent
{\bf Question.} Let $G$ be a finite group embedded in the mapping class group $\Gamma_g$ of genus $g$ and $g\geq 2$. For which pairs $(g,G)$ are the groups 
$\Gamma_{g,G}=\Norm_{\Gamma_g} (G)$ {\it good} groups in the sense of Serre or equivalently by the theorem of Artin-Mazur \cite{AM}, 6.6,
when is the Eilenberg-MacLane space $K(\Gamma_{g, G}^{\wedge}, 1)$ weakly equivalent to the profinite completion $K(\Gamma_{g,G}, 1)^{\wedge}$?

\ms\noindent
In Oda \cite{O} a similar question, attributed to Deligne and Morava was discussed for the mapping class groups $\Gamma_{g, n}$. Recently, Boggi \cite{B}
gave a positive answer to this question for all $\Gamma_{g, n}$ with $2g+n-2>0$, i.e all the mapping class groups are in fact good groups.

\ms\noindent
Using the Grothendieck short exact sequences of etale fundamental groups as in \cite{Z}, Corollary 6.6. (cf. also \cite{O}, \cite{Ma})
we can finally derive the following corollary of our main theorem relating the etale fundamental groups of the moduli stacks $\m_{g,G}$,
the groups $\Gamma_{g,G}$ and the absolute Galois group $\Gal(\bar{\Q}/\Q)$.

\ms
\begin{corollary}
Let $G$ be a finite group embedded in the mapping class group $\Gamma_g$ of genus $g\geq 2$. Let $x$ be any point in the
stack ${\mathcal M}_{g,G}\otimes\bar{\Q}$, then there is a short exact sequence of profinite groups
$$1\rightarrow \Gamma_{g,G}^{\wedge} \rightarrow \pi^{et}_{1}(\m_{g,G}, x) \rightarrow \Gal(\bar{\Q}/\Q)\rightarrow 1.$$
\end{corollary}

\ms\noindent
It also follows from the main theorem that all higher etale homotopy groups of the stacks $\m_{g,G}\otimes \bar{\Q}$ in the sense of Artin-Mazur \cite{AM} are trivial, 
i.e., we have $\pi^{et}_{n}(\m_{g,G}\otimes\bar{\Q}, x)={0}$ for all $n\geq 2$.

\end{document}